# Advanced geometrical constructs in a Pueblo ceremonial site, c 1200 CE

S.Towers[a]

[a]*Simon A. Levin Mathematical, Computational, and Modeling Sciences Center, Arizona State University, Tempe, AZ, USA*

Abstract

Summer 2015 marked the 100th anniversary of the excavation by J.W. Fewkes of the Sun Temple in Mesa Verde National Park, Colorado; an ancient complex prominently located atop a mesa, constructed by the ancestral Pueblo peoples approximately 800 years ago. While the D-shaped structure is generally recognized by modern Pueblo peoples as a ceremonial complex, the exact uses of the site are unknown, although the site has been shown to have key solar and lunar alignments.

In this study, we examined the potential that the site was laid out using advanced knowledge of geometrical constructs. Using aerial imagery in conjunction with ground measurements, we performed a survey of key features of the site. We found apparent evidence that the ancestral Pueblo peoples laid out the site using the Golden rectangle, Pythagorean 3:4:5 triangles, equilateral triangles, and 45 degree right triangles.

The survey also revealed that a single unit of measurement, $L = 30.5 \pm 0.5$ cm, or one third of that, appeared to be associated with many key features of the site. Further study is needed to determine if this unit of measurement is common to other ancestral Pueblo sites, and also if geometric constructs are apparent at other sites.

These findings represent the first potential quantitative evidence of knowledge of advanced geometrical constructs in a prehistoric North American society, which is particularly remarkable given that the ancestral Pueblo peoples had no written language or number system.

[a] smtowers@asu.edu

# Advanced geometrical constructs in a Pueblo ceremonial site, c 1200 CE

GraphicalAbstract

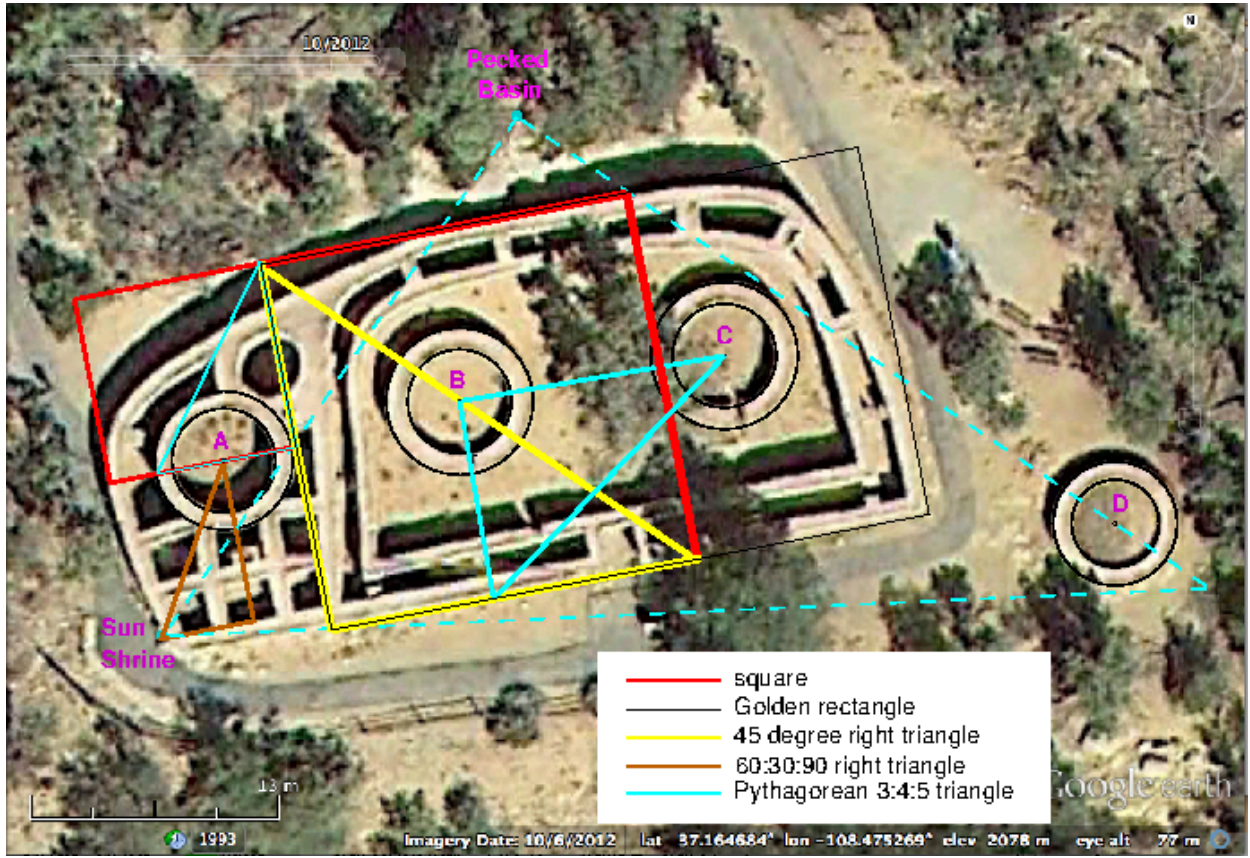

1. Introduction

Mesa Verde National Park is a national park and World Heritage Site, located in southwestern Colorado. The park covers an area of over 210 km$^2$, with the topography consisting of a series of many small mesas separated by deep side canyons [1]. The area was settled by ancestral Pueblo peoples beginning around 470 CE, with final abandonment in late 1200's due to drought conditions [2–6]. Several thousand ruins associated with this period of occupation are found throughout the park [7], with the most famous structure being the Cliff Palace, which was built into a cliff underneath a large rock overhang. The Cliff Palace site was first inhabited in the mid-1000's CE, and finally abandoned in the late 1200's as the region depopulated [8, 9].

Directly across the canyon from Cliff Palace is a site known as the Sun Temple, built atop a mesa with a commanding view of the surrounding landscape. An aerial view of the Sun Temple ruin is shown in Figure 1. The D-shape of the complex is recognized by modern Pueblo peoples to denote a ceremonial structure, however information regarding the exact use of such structures has been lost in oral traditions. The complete lack of domestic artifacts and trash mounds associated with the site point to its use for ceremony, rather than habitation, and the site is extraordinarily unique in the region in this respect, and also in its architecture [10, 11].

The site has proven difficult to date, largely because of the lack of artifacts and wood for dendrochronological dating; however, based upon its masonry style and geographic proximity to the Cliff Palace, the site has been presumed to have been constructed in the 1200's [7, 10, 12].

In addition to the D-shaped outer walls, a notable feature of the Sun Temple complex is the incorporation of four large walled circular structures, which, following the original excavator of the site, Jesse Walter Fewkes [10], we shall refer to as "kivas" (even though, apart from their circular walls, they deviate in many respects from the usual form of a ceremonial kiva in traditional Pueblo architecture [12, 13], and in many respects their original form was more akin to a typical tower [12]). These structures are indicated on the aerial view of the Sun Temple in Figure 1. Following Fewkes, we refer to these structures as Kivas A, B, C, and D [10], labelled in order from west to east; Kiva A is located in the multi-chambered western portion of the complex known as the "Annex", which is connected to the large D-shaped structure. Kivas B and C lie within the D-shaped structure, and Kiva D lies outside the walls to the east. There is another small circular room ("Room s" as indicated in Figure 1) that is markedly different than the four Kivas in several respects; the masonry in the interior of the room is much coarser than in any other area of the site, and it is not constructed of concentric core-and-veneer walls like the four kiva-like structures.

Based on the masonry patterns, the four Kivas have been posited to have preceded the construction of the remainder of the complex, but by an unknown period of time [12]. The construction of the D shaped walls is also posited to have preceded the construction of the walls of the Annex [12].

There is a notable ground feature on the southwest corner, between two short knee-walls that jut out from the side of the complex to either side of a naturally eroded star-shaped basin less than half a meter across. Fewkes dubbed this feature the "Sun Shrine". In addition, a small pecked basin, too small to be visible from aerial imagery (and also surrounded by trees), lies just under 5 meters from the north wall. Two smaller pecked features are located within 30 cm of the pecked basin. The pecked basin and Sun Shrine are indicated in Figure 1.

When Jesse Walter Fewkes excavated the Sun Temple site in the early 1900's, he noted that the walls were made of fine, carefully pecked masonry blocks, were exceptionally vertical, and that their original height was likely around two meters above the present height [10]. He also noted that the complex had had no roof, and erroneously concluded that the structure had never been completed [10, 12]. The extraordinary care with which the site was constructed, and its various unusual architectural features, all point to it being a focus of ceremony in the region.

Fewkes repaired parts of the complex (for instance, by installing capping on the walls to prevent erosion), but from a recent survey of the site and the photographs of the excavation process, it has been concluded that Fewkes did not change the layout of the site [12]. The width of ruin at widest is around 20 m, and length is around 37 meters. The walls are on average 1.3 meters thick, made of masonry surrounding rubble core. The remains of the walls currently range from approximately a meter high for Kivas B, C, and D, to approximately two meters high around the outer D, to approximately two to four meters high for Kiva A and other walls in the Annex.

Previous studies have shown that key architectural features in the Cliff Palace have solar solstice and lunar standstill alignments with the Sun Temple [12, 14], including the Kivas and the Sun Shrine. While using aerial imagery to perform site measurements to examine these alignments, the author incidentally noted several apparent repeated measures, and some remarkable apparent geometrical constructs in the site, which led to this work; here we present a detailed analysis of apparent geometrical constructs at the site related to the positions of the outer walls enclosing the site, the location of the Sun Shrine, and the radii and position of the centers of the four Kivas. The analysis is based on both the ground and aerial surveys, because the metrology of some site features is difficult to precisely and accurately assess on the ground due to the uneven nature of the site walls. Thus, not only does the ground survey serve as a cross-check of the accuracy of the aerial survey, but the measurements of the two surveys also complement each other.

With a site survey based on both aerial imagery and ground measurements of the site, we examined the possibility that squares, Pythagorean 3:4:5 right triangles, equilateral triangles (or 30°:60°:90° right triangles), and Golden rectangles were incorporated into the layout of the site. In the analysis, we also examined the potential that a common unit of measure underlies the observed apparent geometric shapes.

There have been few prior comprehensive studies of geometrical constructs in architecture in the prehistoric New World, and have largely been confined to examination of Mayan architecture [15–17]. Thus, this study will help shed needed light on the scope of geometrical knowledge in a prehistoric North American society. In addition, to our knowledge, no analysis has previously been published that examines the potential of a common unit of measurement underlying the layout of any prehistoric site in the Southwest.

In the following sections, we describe the aerial and ground surveys of the site. We also provide a brief description of the geometric properties of Pythagorean 3:4:5 right triangles, equilateral triangles, and Golden rectangles, followed by a presentation and discussion of results.

*Figure 1:* Aerial view of the Sun Temple complex (as obtained from Google Earth, accessed May 1, 2016), with pertinent features labeled according to [10], and the location of the Sun Shrine indicated (a naturally eroded star-shaped basin, bracketed by two small knee-walls on the SW corner of the site). Also indicated is the approximate position of the small pecked basin to the north of the site (too small to be visible in aerial imagery, and also surrounded by trees). The position of the pecked basin is estimated from ground survey measurements. The ground length over the horizontal width of the view is approximately 65 meters. Google and the Google logo are registered trademarks of Google Inc., used with permission.

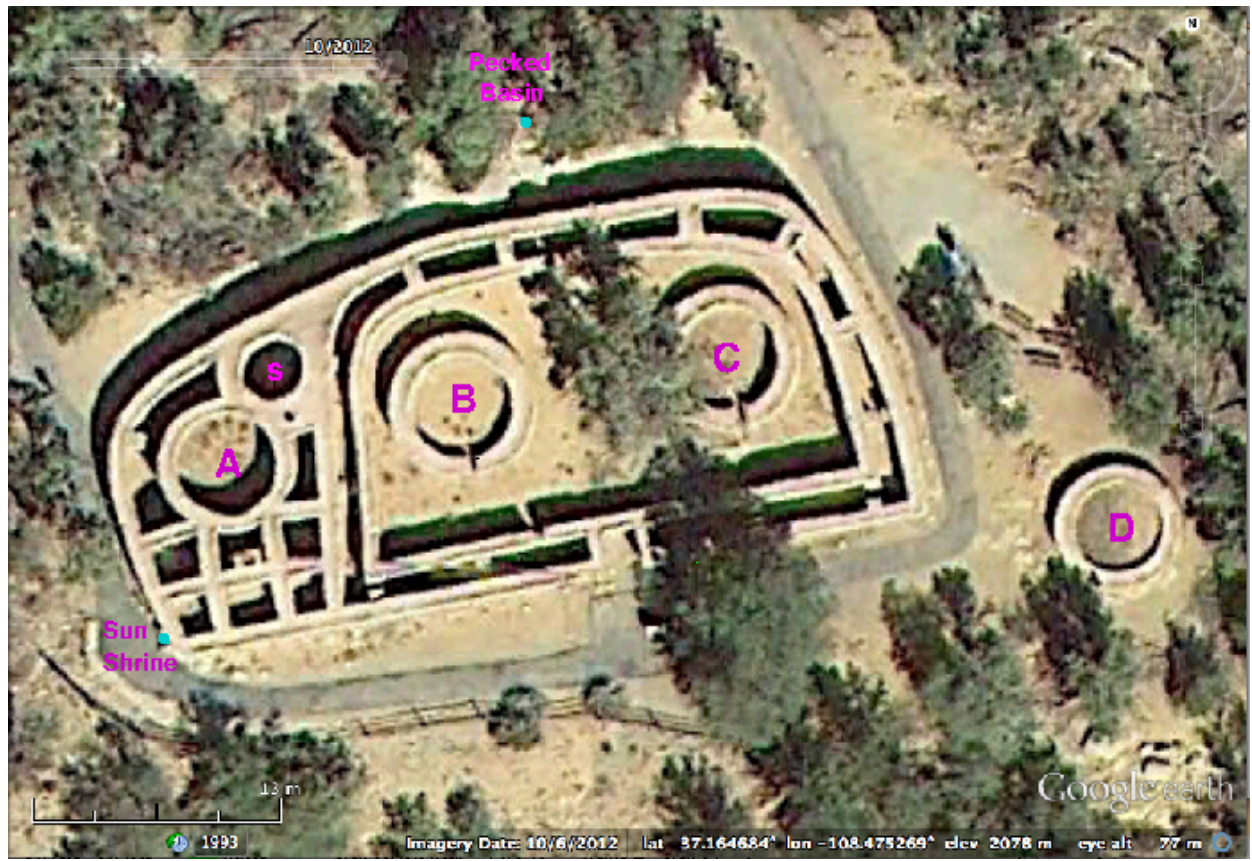

2. Methods and Materials

*2.1. Examination of site for potential geometrical constructs, and uncertainties in measurements*

Using aerial and ground surveys, described below, we examined the dimensions of the rectangle that encased the outer D of the site, and the dimensions of the radii of the four Kivas. We also examined the relative lengths and angles of lines between key site features, to determine if there was any consistency with 90°, 45°, or 60° angles, or the angles associated with the Pythagorean 3:4:5 triangle. We also examined the potential that a common unit of measurement underlay the site layout.

We constrained the analysis to examination of site elements that were found by previous analyses to have been built first [12]. This includes the four Kivas, the outer D, the Sun Shrine, and the pecked basin. We note here that all the walls of the Sun Temple complex are exceptionally vertical, but the masonry of the exterior of the wall of the outer D is somewhat finer and more straight than that of the interior side of the wall, indicating more precise attention to detail in its construction than that seen in the interior wall. In addition, masonry patterns at the SW and SE corners on the exterior surface of the outer D indicate that the builders took pains to ensure that the walls met at an almost exact 90° angle at those corners, but this is not true of the interior surface of the wall. For these reasons, we constrained the analysis to points on the exterior wall only, and did not consider measurements relative to the interior side of the wall. For the outer D, we thus only considered measurements related to the length and width of the rectangle that encases it, and the positions of other site features relative to the SE and SW corner. Other geometrical constructs may be associated with other site features that we did not consider in this analysis; further study may or may not reveal such geometries.

There were two sources of uncertainties associated with the measurements recorded in this analysis. The first were the measurement uncertainties associated with our surveys due to issues such as the uneven nature of the site walls, aerial image pixelation, tape measure stretch or sag, and/or theodolite uncertainty measurements. Repeated measures of each feature were used in both surveys to estimate the uncertainty associated with the measurement. The second were the uncertainties associated with the measurements of the original site layout; if the ancestral Pueblo peoples used a common unit of measure to construct the site, there was some uncertainty associated with replication of this measure over and over, particularly if the span of a body part (like a foot or hand) was used, or a cord was used for some measurements, because cords can stretch. Even with careful attention to detail, there would have been some measurement uncertainty associated with the site layout.

### 2.2. Digital survey using aerial imagery

Google Earth is a virtual globe, map, and geographical information systems (GIS) program, freely available from `http://earth.google.com` (accessed December 31, 2016). Since the launch of the product in 2005, it has been used in a wide range of academic endeavors, including for use in the survey of archaeological sites (see, for instance, References [18–22]). In this analysis, we used Google Earth to obtain aerial imagery and geographic information related to the Sun Temple site.

To survey the site, an aerial view of the site was obtained from Google Earth, including the image distance scale (see, for example, Figure 1). Google Earth allows selection through all past aerial images, enabling the user to select an image with the best contrast of site features. The image was then read into Xfig, a free and open-source vector graphics software package[b].

Within Xfig, circles were overlaid onto the four Kivas, and the radii and centers of the circles in the coordinate frame of the image were determined, as were the length of lines connecting various key site features. Because there was some amount of objectivity involved in the placement of the circles and lines, the procedure was repeated several times, and the average and one standard deviation uncertainty on the measurements determined from the iterations.

### 2.3. Ground Survey

Upon obtaining a research permit from the National Park Service, the author visited the Sun Temple site in summer 2015, and performed a survey of various features of the site using theodolite and tape measurements to verify the aerial survey measurements (also known as "ground truthing" [23]). Measurements of the dimensions of the outer D, and of the inner and outer diameters of the four Kivas were obtained, as was the distance between Kivas B and C, and the distances of Kivas B, C, and D from the outer D wall enclosing the site complex.

The walls of Kivas B, C, and D are approximately a meter or less high, whereas the walls in the vicinity of Kiva A are substantially higher (between three to four meters). The walls of Kiva A slope slightly inwards, thus the radius of Kiva A is slightly larger at the ground level than at the top, by approximately 5 cm.

---

[b] See `www.xfig.org`, accessed December 31, 2016

Due to the uneven nature of the walls, issues of tape stretch and/or sag over long distances, theodolite precision, and the somewhat deteriorated nature of some of the walls on the site, the ground measurements have some uncertainty which must be taken into account in statistical comparisons of the ground measurements to the aerial survey measurements.

*2.4. Geometrical constructs: Pythagorean and Equilateral triangles, and the Golden rectangle*

In this analysis, we examined the layout of the Sun Temple site for evidence of squares, 45° right triangles, Pythagorean 3:4:5 triangles, equilateral triangles (or 60°:30°:90° right triangles), and Golden rectangles.

*Pythagorean Triangles*

Pythagorean triples are integers x, y, and z, such that $z^2 = x^2 + y^2$. A right triangle whose sides are Pythagorean triples is known as a Pythagorean triangle, and such triangles were known to several ancient societies in Asia, the Middle East, and the Mediterranean [24]. There are infinitely many Pythagorean triples, with the lowest order triple being (x, y, z) = (3, 4, 5). Use of the 3:4:5 triangle is a standard practice in construction in modern times as a simple means to obtain walls at right angles.

*Equilateral triangles*

Equilateral triangles have equal length on all three sides, and interior angles all equal to 60°. Using a straightedge and a compass, it is straightforward to construct equilateral triangles of a given side length, as shown in Figure 2. Equilateral triangles of a given height can also be readily constructed [25].

Right triangles with angles 60°:30°:90° are easily obtained by halving an equilateral triangle, as shown in Figure 2. The ratios of the side lengths of such triangles are 1:√3:2.

*Golden rectangle*

The Golden ratio (or Golden mean, or Golden section) is

$$\phi = \frac{(1 + \sqrt{5})}{2} \sim 1.618$$

and is a ratio found often in nature, and also employed in ancient and modern architecture, as rectangles with ratio of side lengths equal to $\phi$ are felt to be 'pleasing' in appearance [26]. As shown in Figure 2, Golden rectangles are readily constructed with a straightedge and a compass [25], and do not require knowledge of irrational numbers.

*Figure 2:* The upper left plot shows the construction of an equilateral triangle, beginning with the green side, of side length 1. Circles of radius 1 are drawn with a cord from each end of the green line. Their intersection forms the apex of the triangle. A 60°:30°:90° right triangle is obtained by drawing a line from the apex to the point midway on the green line. The upper right plot shows the details of the construction of a Golden rectangle of width 1. The lower left plot shows the dimensions of two circles inscribed and circumscribed on a square with side length equal to 2. The walls of Kivas B, C, and D appear to be constructed this manner (see Table 2).

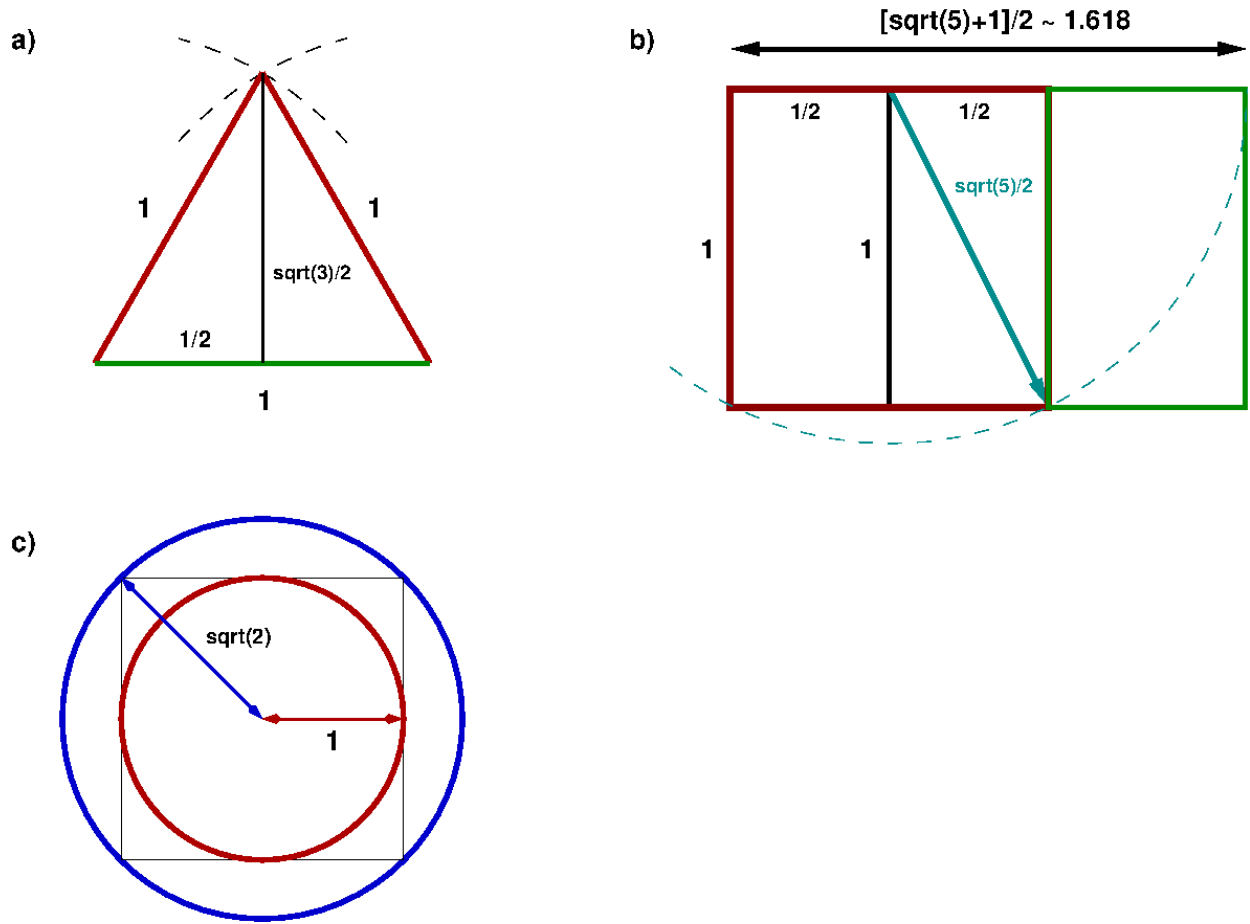

3. Results

The aerial and ground survey measurements of the features are summarized in Table 1. In the following subsections, we describe the specifics of each of the measurements, and how they related to one and other, either as evidence of apparent geometrical constructs, or evidence for an apparent common unit of measurement.

Prior to those descriptions, however, it must be noted that many statistical hypothesis tests were performed in this analysis, and this has implications for tests of significance; for instance, there are 17 tests of significance comparing the aerial and ground survey measurements in Table 1 alone. Because, by mere random chance, the p-value testing the null hypothesis will be $p < 0.05$ five percent of the time when the null hypothesis is actually true, the cut-off, $\alpha = 0.05$, for rejecting the null hypothesis must be adjusted when multiple tests of significance are performed [27]. The Bonferroni correction uses the cut-off $\alpha' = \alpha/k$, where k is the number of tests of significance [27]. In this analysis, we performed almost four dozen tests of significance, thus we rejected the null hypothesis only when $p < \sim 0.001$.

In all cases examined, the aerial and ground survey measurements of site features are statistically consistent.

*3.1. Apparent geometrical constructs*

In Figure 3, we show the Sun Temple site, with several apparent geometrical constructs overlaid.

*3.1.1. Golden rectangle, and associated squares*

The aerial and ground survey measurements of the rectangle encasing the outer D are shown in lines 9 and 10 of Table 1. We estimate from the aerial and ground surveys that the ratios of the length to width of the rectangle encasing the outer D are $1.643 \pm 0.013$, and $1.646 \pm 0.013$, respectively, where the uncertainties shown are the one standard deviation uncertainties. These ratios of the length to width of the outer D are both statistically consistent with the Golden ratio, $\phi \sim 1.618$ ($\chi^2_{\nu=1}$ p=0.06 and p=0.03, respectively).

As seen in Figure 3, the west side of the rectangle encasing the outer D appears to be tangent to the outer radius of Kiva A, and the line perpendicular to the midpoint of the rectangle goes through the center of Kiva A. If X is the width of the rectangle encasing the outer D of the site (i.e. the length of the red lines in Figure 4), the line through the center of Kiva A that bisects X at a perpendicular, when extended a distance of X/2 to form the bottom of a square, touches the furthest western edge of the Annex walls.

*3.1.2. Equilateral triangle*

As seen in Figure 3, the line between the Sun Shrine and the center of Kiva A makes an angle of 60° with the line along the inner wall of the base of the D, to within less than a degree. This may be evidence of an equilateral triangle, or, alternatively, it may be the hypotenuse of a 1:√3:2 right triangle (i.e. a triangle with interior angles 60°:30°:90°, as indicated in Figure 3). Indeed, such a triangle has one side that goes through the ventilator structure immediately to the south of Kiva A.

*3.1.3. Pythagorean 3:4:5 triangles*

As shown in Figure 3, there is an apparent Pythagorean 3:4:5 triangle associated with the construction of the outer radius of Kiva A relative to the side of the rectangle that encases the outer D. From the aerial and ground survey measurements, shown in lines 1 and 10 of Table 1, we estimate that the ratios of the width of the outer D to the outer radius of Kiva A are $5.337 \pm 0.060$ and $5.397 \pm 0.061$, respectively. These values are both close to $16/3 \sim 5.333$, as they would be if constructed from the side of the westernmost Pythagorean 3:4:5 triangle shown in Figure 3. The p-values testing the null hypothesis that these are statistically consistent with $16/3 = 5.333$ are $p = 0.95$ and $p = 0.30$, respectively ($\chi^2_{\nu=1}$).

The ratio of the outer to inner radius of Kiva A is statistically consistent to 4/3, as shown in Table 2. This may indicate evidence of an additional Pythagorean 3:4:5 triangle used in the construction of the walls of Kiva A.

There is an apparent Pythagorean 3:4:5 triangle associated with the centers of Kivas B and C, and the nearest distance of the center of Kiva B to the south wall of the outer D. The aerial and ground survey measurements of the distance between centers of Kivas B and C, and the distance of the center Kiva B to south wall of outer D, perpendicular to the line between the centers of Kivas B and C, are shown in lines 12 and 13 of Table 1. From the aerial and ground survey measurements, the estimates of the ratios of the two measurements are $1.363 \pm 0.015$ and $1.350 \pm 0.015$, respectively.

These two estimates are close to 4/3, as they would be if the two lines are adjacent to the right angle in a Pythagorean 3:4:5 triangle. The p-values testing the null hypothesis that these are statistically consistent with 4/3 are $p = 0.05$ and $p = 0.28$, respectively ($\chi^2_{\nu=1}$).

### 3.1.4. Circles circumscribed on squares, and vice versa

The aerial and ground survey measurements of the inner and outer radii of the four Kivas are seen in lines 1 to 8 of Table 1.

The estimates of the ratios of the outer to inner radii of the Kivas from the aerial and ground surveys are shown in Table 2. For the ratio of the outer to inner radii of Kiva A, we added 5 cm to the radii obtained from examination of the top of the 4 meter high structure, because ground measurements reveal the walls slope slightly inwards as they rise.

The ratios of the outer to inner radii of Kivas B, C, and D are close to $\sqrt{2}$, which would be obtained if the walls were constructed from circles inscribed and circumscribed on a square, as shown in Figure 2. In Table 2 we compare the ratios of the radii of Kivas B, C, and D to $\sqrt{2}$. In all cases, the measurements are statistically consistent with these quantities.

### 3.2. Other consistencies in site layout

The line going from the Sun Shrine parallel to the south wall of the outer D is tangent to the outer wall of Kiva D.

In addition, the center of Kiva B lies on one of the 45° diagonals of the square used in the construction of the Golden rectangle encasing the outer D.

As seen in Figure 3, the diagonal that goes through Kiva B also appears to be perpendicular, to within less than a degree, to the line that goes between the Sun Shrine and the pecked basin. The line from the Sun Shrine to the pecked basin also appears to be perpendicular, to within less than a degree, to the line that goes from the pecked basin to the center of Kiva D. These two lines may form two sides of a Pythagorean triangle, as shown in Figure 3. Munson *et al* (2010) mention that they believe the pecked basin may have been a site datum point [12], and this is supported by the geometric patterns we observe in this analysis. Indeed, within 30 cm of the pecked basin, there are two smaller pecked features which may have been initial attempts at placement of the datum point, but were not completed when the mistakes were detected; the two smaller features do not appear to be associated with any geometric constructs at the site.

*Table 1: Summary of the aerial and ground survey measurements of the inner and outer radii of the four Kivas, and other key site features. The uncertainties shown are the one standard deviation uncertainties. The last column is the p-value testing the null hypothesis that the true quantities underlying the aerial and ground and survey estimates are equal.*

|    |                                                                  | Aerial imagery survey (cm) | Ground survey (cm) | p-value $\chi^2_{\nu=1}$ |
|----|------------------------------------------------------------------|----------------------------|--------------------|--------------------------|
| 1  | Inner radius Kiva A                                              | 265±2                      | 264±2              | p=0.72                   |
| 2  | Outer radius Kiva A (at top; +5cm for ground radii)              | 360±3                      | 355±3              | p=0.24                   |
| 3  | Inner radius Kiva B                                              | 271±2                      | 267±2              | p=0.16                   |
| 4  | Outer radius Kiva B (at ground)                                  | 385±3                      | 382±3              | p=0.48                   |
| 5  | Inner radius Kiva C                                              | 268±2                      | 265±2              | p=0.29                   |
| 6  | Outer radius Kiva C (at ground)                                  | 382±3                      | 383±3              | p=0.81                   |
| 7  | Inner radius Kiva D                                              | 233±2                      | 235±2              | p=0.48                   |
| 8  | Outer radius Kiva D (at ground)                                  | 332±3                      | 334±3              | p=0.64                   |
| 9  | Distance along south wall of outer D                             | 3200±8                     | 3199±8             | p=0.93                   |
| 10 | Width of outer D                                                 | 1948±15                    | 1943±15            | p=0.81                   |
| 11 | Nearest distance between Kivas B and C                           | 653±8                      | 643±5              | p=0.29                   |
| 12 | Distance between centers Kivas B and C                           | 1422±8                     | 1417±8             | p=0.66                   |
| 13 | Distance center of Kiva B approx. $\perp$ to south all of outer D | 1043±10                   | 1050±10            | p=0.62                   |
| 14 | Nearest distance outer radius of Kiva B to SW corner of outer D  | 974±4                      | 972±8              | p=0.82                   |
| 15 | Nearest distance outer radius of Kiva C to SE corner of outer D  | 971±5                      | 960±8              | p=0.24                   |
| 16 | Nearest distance outer radius of Kiva D to SE corner of outer D  | 671±8                      | 655±7              | p=0.13                   |
| 17 | Distance center of Kiva D to SE corner of outer D                | 984±5                      | 998±8              | p=0.14                   |
| 18 | Distance Sun Shrine to center of Kiva A                          | 972±10                     | -                  | -                        |
| 19 | Distance center of Kiva A to south wall                          | 967±10                     | -                  | -                        |

Accepted *Journal of Archaeological Science: Reports*, Jan 2017 (in press)

*Table 2: Ratio of the outer to inner ground radii of the four Kivas, based on the aerial and ground survey measurements, and compared to the hypothesized values.*

|        | Aerial survey ratio outer to inner radii | Comparison | p-value $\chi^2_{v=1}$ |
|--------|------------------------------------------|------------|------------------------|
| Kiva A | 1.352±0.015                              | 4/3~1.333  | p=0.22                 |
| Kiva B | 1.421±0.015                              | √2~1.414   | p=0.67                 |
| Kiva C | 1.425±0.015                              | √2~1.414   | p=0.47                 |
| Kiva D | 1.425±0.018                              | √2~1.414   | p=0.55                 |

|        | Ground survey ratio outer to inner radii | Comparison | p-value $\chi^2_{v=1}$ |
|--------|------------------------------------------|------------|------------------------|
| Kiva A | 1.338±0.015                              | 4/3~1.333  | p=0.74                 |
| Kiva B | 1.431±0.016                              | √2~1.414   | p=0.29                 |
| Kiva C | 1.445±0.016                              | √2~1.414   | p=0.05                 |
| Kiva D | 1.421±0.018                              | √2~1.414   | p=0.69                 |

### 3.3. Measurements of key features of the site, and evidence for a common unit of measurement

In Figure 4, we show an aerial view of the Sun Temple site, with lines overlaid representing measurements of several key features of the site. In the figure, the yellow lines are set to exactly one half the length of the red lines (which represent the width of the rectangle encasing the outer D), and the dark blue lines are set to exactly one third the length of the red lines. The brown line is set to exactly 3/8 the length of the red, and the pink circles overlaying the inner radii of Kivas A, B and C have exactly the same radii. Repeated measures of a common unit appear to be evident in the relative positions of many of the key site features. In the following sections, we discuss each of these in turn.

#### 3.3.1. Consistency of inner radii of Kivas A, B, and C

The ground and aerial survey measurements of the inner radii of Kivas A, B, and C are shown in lines 1 to 6 of Table 1. In order to estimate the inner radius of Kiva A at ground level, we added 5 cm to the radii obtained from examination of the top of the ~4 meter high structure, because ground measurements reveal the walls of Kiva A slope slightly inwards as they rise. The remnants of the walls of Kivas B and C are significantly shorter and close to ground level.

The aerial and ground survey measurements of the inner radii of Kivas A, B, and C at the ground level are statistically consistent with being equal for the three Kivas, with mean $270 \pm 2$ cm and $267 \pm 2$ cm, respectively ($\chi^2_{\nu=1}$ p-values $p = 0.56$ and $p = 0.37$, respectively).

#### 3.3.2. Ratio of the width of rectangle encasing the outer D to the outer radius of Kiva A

From the aerial and ground survey measurements, shown in lines 2 and 10 of Table 1, we estimate that the ratios of the width of the outer D to the outer radius of Kiva A are $5.337 \pm 0.060$ and $5.397 \pm 0.061$, respectively. These values are both close to 16/3 ~ 5.333, as they would be if constructed from the side of the westernmost Pythagorean 3:4:5 triangle shown in Figure 3. The p-values testing the null hypothesis that these are statistically consistent with $16/3 = 5.333$ are $p = 0.95$ and $p = 0.30$, respectively ($\chi^2_{\nu=1}$).

*3.3.3. Ratio of the width of rectangle encasing outer D to the inner radius of Kiva A*

From the aerial and ground survey measurements, shown in lines 1 and 10 of Table 1, we estimate that the ratios of the width of the outer D to the inner radius of Kiva A are $7.215 \pm 0.077$ and $7.223 \pm 0.077$, respectively. These values are both close to $64/9 \sim 7.111$, as they would be if the outer radius was $4/3$ the inner radius, and the outer radius was constructed from the side of the westernmost Pythagorean 3:4:5 triangle shown in Figure 3. The p-values testing the null hypothesis that these are statistically consistent with $64/9 = 7.111$ are $p = 0.18$ and $p = 0.15$, respectively ($\chi^2_{v=1}$).

*3.3.4. Ratio of the width of the rectangle encasing outer D to the outer radius of Kiva D*

From the aerial and ground survey measurements, shown in lines 8 and 10 of Table 1, we estimate that the ratios of the width of the outer D to the outer radius of Kiva D are $5.867 \pm 0.070$ and $5.817 \pm 0.069$, respectively. These two values are both close to six, perhaps indicating a common underlying unit of measure. The p-values testing the null hypothesis that these are statistically consistent with 6 are $p = 0.06$ and $p = 0.01$, respectively ($\chi^2_{v=1}$).

*3.3.5. Ratio of the width of the rectangle encasing the outer D to the nearest distance between the outer radius of Kiva B and the SW corner of the outer D*

From the aerial and ground survey measurements, shown in lines 10 and 14 of Table 1, we estimate that the ratios of the width of the outer D to the nearest distance between the outer wall of Kiva B to the SW corner of the outer D are $2.000 \pm 0.017$ and $1.999 \pm 0.023$, respectively. These two values are both close to two, perhaps indicating a common underlying unit of measure. The p-values testing the null hypothesis that these are statistically consistent with 2 are $p = 1.00$ and $p = 0.96$, respectively ($\chi^2_{v=1}$).

*3.3.6. Ratio of the width of the rectangle encasing the outer D to the nearest distance between the outer radius of Kiva C and the SE corner of the outer D*

From the aerial and ground survey measurements, shown in lines 10 and 15 of Table 1, we estimate that the ratios of the width of the outer D to the nearest distance between the outer wall of Kiva C to the SE corner of the outer D are $2.006 \pm 0.019$ and 2.024 ± 0.023, respectively. These two values are both close to two, perhaps indicating a common underlying unit of measure. The p-values testing the null hypothesis that these are statistically consistent with 2 are p = 0.74 and p = 0.30, respectively ($\chi^2_{\nu=1}$).

*3.3.7. Ratio of the width of the rectangle encasing the outer D to the nearest distance between the outer radius of Kiva D and the SE corner of the outer D*

From the aerial and ground survey measurements, shown in lines 10 and 16 of Table 1, we estimate that the ratios of the width of the outer D to the nearest distance between the outer wall of Kiva D to the SE corner of the outer D are 2.903 ± 0.041 and 2.966 ± 0.039, respectively. These two values are both close to three, perhaps indicating a common underlying unit of measure. The p-values testing the null hypothesis that these are statistically consistent with 3 are p = 0.02 and p = 0.39, respectively ($\chi^2_{\nu=1}$).

*3.3.8. Ratio of the width of the rectangle encasing the outer D to the distance between the center of Kiva D and the SE corner of the outer D*

From the aerial and ground survey measurements, shown in lines 10 and 16 of Table 1, we estimate that the ratios of the width of the outer D to the distance of the center of Kiva D to the SE corner of the outer D are 1.980 ± 0.018 and 1.947 ± 0.022, respectively. These two values are both close to two, perhaps indicating a common underlying unit of measure. The p-values testing the null hypothesis that these are statistically consistent with 2 are p = 0.27 and p = 0.01, respectively ($\chi^2_{\nu=1}$).

*3.3.9. Ratio of the width of the rectangle encasing the outer D to the nearest distance between the outer radius of Kiva B and the outer south wall of the outer D*

From the aerial and ground survey measurements, shown in lines 10, 4, and 13 of Table 1, we estimate that the ratios of the width of the outer D to the nearest distance between the outer wall of Kiva B to the south wall of the outer D are 2.960 ± 0.050 and 2.909 ± 0.049, respectively. These two values are both close to three, perhaps indicating a common underlying unit of measure. The p-values testing the null hypothesis that these are statistically consistent with 3 are p = 0.42 and p = 0.06, respectively ($\chi^2_{\nu=1}$).

*3.3.10. Ratio of the width of the rectangle encasing the outer D to the nearest distance between the outer radii of Kivas B and C*

From the aerial and ground survey measurements, shown in lines 10 and 11 of Table 1, we estimate that the ratios of the width of the outer D to the nearest distance between the outer radii of Kivas B and C are $2.983 \pm 0.043$ and $3.022 \pm 0.033$, respectively. These two values are both close to three, perhaps indicating a common underlying unit of measure. The p-values testing the null hypothesis that these are statistically consistent with 3 are $p = 0.70$ and p=0.51, respectively ($\chi^2_{\nu=1}$).

### 3.3.11. Ratio of the width of rectangle encasing outer D to the distance between the Sun Shrine and the center of Kiva A

From the aerial survey measurements, shown in lines 10 and 18 of Table 1, we estimate that the ratio of the width of the outer D to the distance from the Sun Shrine to the center of Kiva A is $2.004 \pm 0.026$. This is close to two, perhaps indicating a common underlying unit of measure. The p-value testing the null hypothesis that this is statistically consistent with 2 is $p = 0.87$ ($\chi^2_{\nu=1}$).

### 3.3.12. Ratio of the width of rectangle encasing outer D to the distance between the center of Kiva A and the south wall

From the aerial survey measurements, shown in lines 10 and 19 of Table 1, we estimate that the ratio of the width of the outer D to the distance from the Sun Shrine to the center of Kiva A is $2.014 \pm 0.026$. This is close to two, perhaps indicating a common underlying unit of measure. The p-value testing the null hypothesis that this is statistically consistent with 2 is $p = 0.59$ ($\chi^2_{\nu=1}$).

*Figure 3: Some apparent geometrical constructs evident in the layout of the Sun Temple site, including squares, 45° right triangles, a 60°:30°:90° right triangle (or an equilateral triangle), a Golden rectangle, and Pythagorean 3:4:5 triangles. Google and the Google logo are registered trademarks of Google Inc., used with permission.*

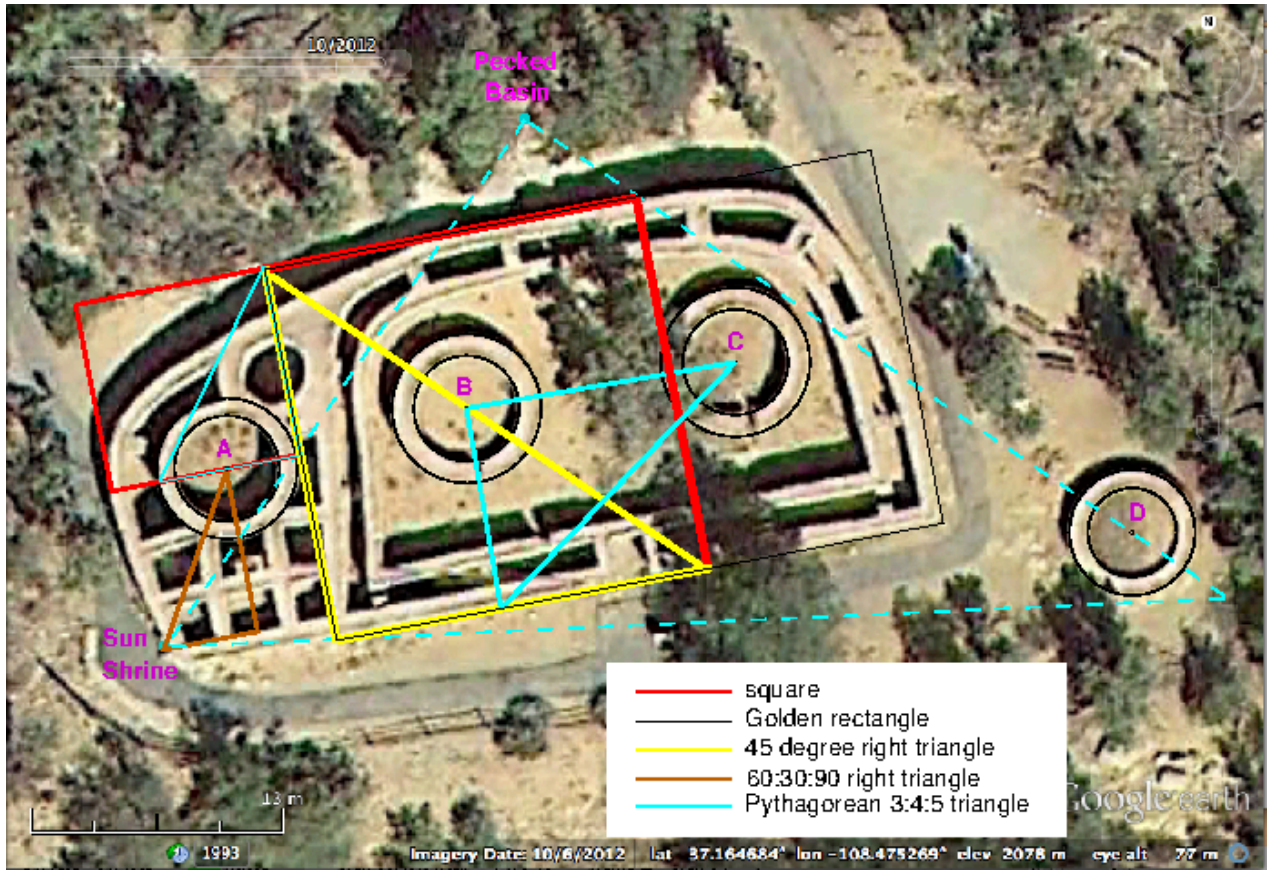

### 3.4. Estimation of the base unit of measurement used to construct the site

In the following, we refer to X as the width of the rectangle encasing the outer D of the site (i.e. the length of the red lines in Figure 4). Based on the observations of apparent consistencies in site measurements and geometrical constructs, we can estimate the value of X and its uncertainty by taking the average and standard deviation of the following:

- The width of the rectangle encasing the outer D.
- $1/\phi$ times the length of the rectangle encasing the outer D
- 64/9 times the inner radius of Kiva A
- 64/9 times the inner radius of Kiva B
- 64/9 times the inner radius of Kiva C
- $6/\sqrt{2}$ times the inner radius of Kiva D
- 16/3 times the outer radius of Kiva A
- $64/(9\sqrt{2})$ times the outer radius of Kiva B
- $64/(9\sqrt{2})$ times the outer radius of Kiva C
- 6 times the outer radius of Kiva D
- 2 times the nearest distance between the outer wall of Kiva B and the SW corner of the outer D
- 2 times the nearest distance between the outer wall of Kiva C and the SE corner of the outer D
- 3 times the nearest distance between the outer wall of Kiva D and the SE corner of the outer D
- 3 times the nearest distance between the outer wall of Kiva B and the south wall of the outer D
- 2 times the distance between the center of Kiva D and the SE corner of the outer D
- 3 times the nearest distance between the two outer walls of Kivas B and C
- 2 times distance between the center of Kiva A and the Sun Shrine
- 2 times distance between the center of Kiva A and the south wall

The results are shown in Tables 3 and 4 for the aerial and ground surveys, respectively. The aerial survey estimates that X = 1952 ± 26 cm, and the ground survey estimates that X = 1945 ± 37 cm. In both cases, the standard deviation uncertainty is less than 2% of X. The variance on the estimate of X represents the variance of our measurements, plus the variance of the measurements of the Pueblo peoples when laying out the site. The average relative standard deviation uncertainty on our measurements is approximately 1%, thus the estimated average relative standard deviation uncertainty on the ancestral Pueblo measurements was also approximately 1% (i.e. an uncertainty of approximately one centimeter in each meter).

Because the inner radii of Kivas A, B, and C are statistically consistent with $9X/64$, it appears that, if these inner radii were laid out as integer increments of some base unit, L, then L is, at most, $L = X/64$. From the aerial and ground surveys, we thus estimate $L = 30.50 \pm 0.41$ cm, and $L = 30.39 \pm 0.58$ cm, respectively. Note that the true base unit measure could either be L, or an integer fraction of L. Indeed, the fact that the nearest distance between the outer walls of Kivas B and C is consistent with $X/3$, as is the nearest distance of the outer radius of Kiva D to the SE corner of the outer D, points to X likely also being divisible by three (thus L is likely three times a sub-unit, equal to approximately 10 cm).

*Figure 4: Measurements related to key features of the Sun Temple site, exhibiting apparent evidence of a common unit of measure. In the figure, the length of all yellow lines is exactly 1/2 that of the red, and the length of all dark blue lines is exactly 1/3 that of the red. The brown line is set to exactly 3/8 the length of the red. The pink circles shown overlaying the inner radii of Kivas A, B and C have the same radii. Google and the Google logo are registered trademarks of Google Inc., used with permission.*

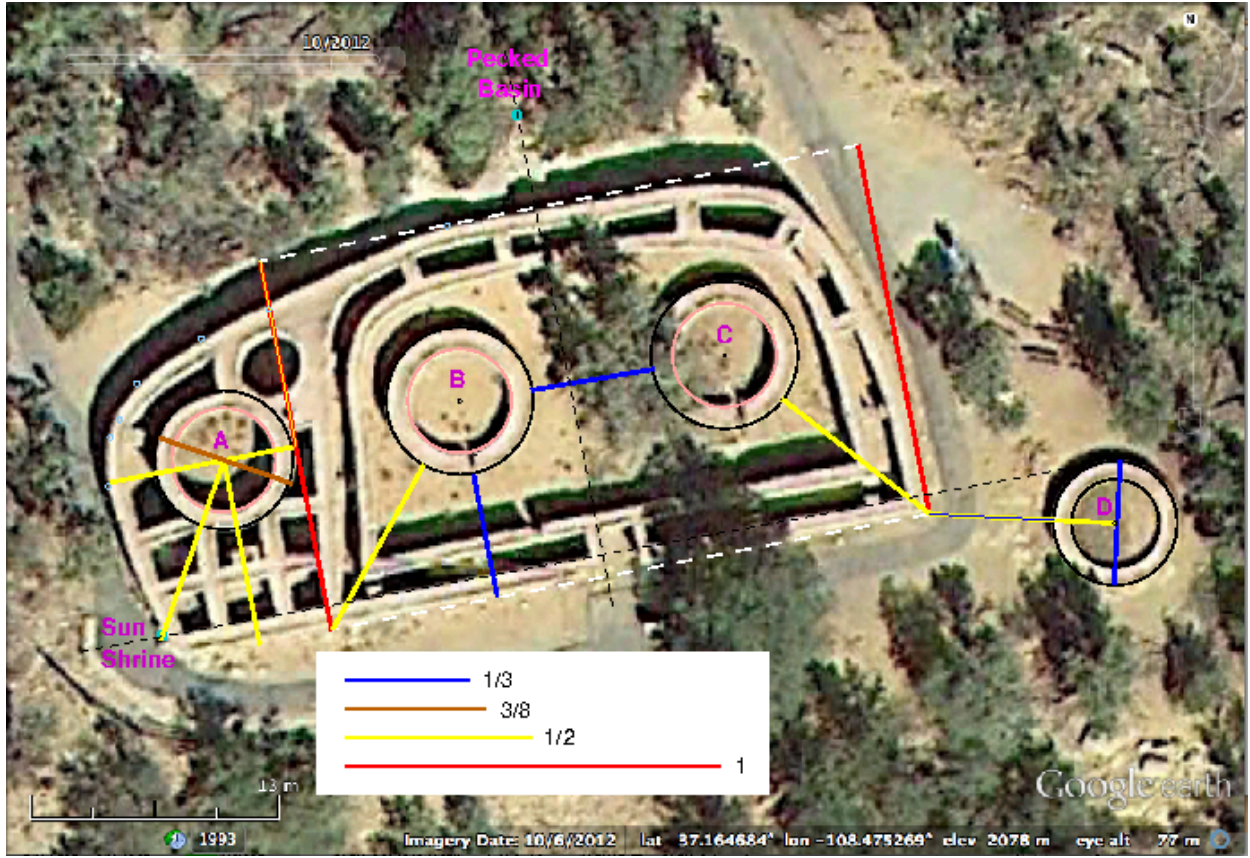

*Table 3: Aerial survey measurements between key features of the site, expressed in the hypothesized relationship to the unit of measure representing the width of the rectangle encasing the outer walls of the D (the multiplying factors are derived from the hypothesized geometrical constructs at the site, seen in Figure 3). The average is weighted according to the uncertainty of the individual measurements.*

|  | Aerial survey (cm) |
|---|---|
| Width of the rectangle encasing outer D | 1948±15 |
| 1/ϕ times the length of the rectangle encasing outer D | 1978±5 |
| 64/9 times the inner radius of Kiva A | 1920±14 |
| 64/9 times the inner radius of Kiva B | 1947±16 |
| 64/9 times the inner radius of Kiva C | 1927±14 |
| 6/√2 times the inner radius of Kiva D | 1936±15 |
| 16/3 times the outer radius of Kiva A | 1906±14 |
| 64/(9/√2) times the outer radius of Kiva B | 1921±15 |
| 64/(9/√2) times the outer radius of Kiva C | 1977±17 |
| 6 times the outer radius of Kiva D | 1992±18 |
| 2 times distance outer wall Kiva B to SW corner of outer D | 1948±8 |
| 2 times distance outer wall Kiva C to SE corner of outer D | 1942±10 |
| 3 times distance outer wall Kiva D to SE corner of outer D | 2013±24 |
| 3 times distance outer wall Kiva B to south wall of outer D | 1974±30 |
| 2 times distance center Kiva D to SE corner of outer D | 1968±10 |
| 3 times nearest distance between outer walls of Kivas B and C | 1959±24 |
| 2 times distance center Kiva A to Sun Shrine | 1944±20 |
| 2 times distance center Kiva A to south wall | 1934±20 |
| Weighted average | 1952±26 |

*Table 4: Ground survey measurements between key features of the site, expressed in the hypothesized relationship to the unit of measure representing the width of the rectangle encasing the outer walls of the D (the multiplying factors are derived from the hypothesized geometrical constructs at the site, seen in Figure 3). The average is weighted according to the uncertainty of the individual measurements.*

|  | Ground survey (cm) |
|---|---|
| Width of the rectangle encasing outer D | 1943±15 |
| 1/φ times the length of the rectangle encasing outer D | 1977±5 |
| 64/9 times the inner radius of Kiva A | 1913±14 |
| 64/9 times the inner radius of Kiva B | 1920±16 |
| 64/9 times the inner radius of Kiva C | 1899±14 |
| 6/√2 times the inner radius of Kiva D | 1921±15 |
| 16/3 times the outer radius of Kiva A | 1884±14 |
| 64/(9/√2) times the outer radius of Kiva B | 1926±15 |
| 64/(9/√2) times the outer radius of Kiva C | 1994±17 |
| 6 times the outer radius of Kiva D | 2004±18 |
| 2 times distance outer wall Kiva B to SW corner of outer D | 1944±16 |
| 2 times distance outer wall Kiva C to SE corner of outer D | 1920±16 |
| 3 times distance outer wall Kiva D to SE corner of outer D | 1965±21 |
| 3 times distance outer wall Kiva B to south wall of outer D | 2004±30 |
| 2 times distance center Kiva D to SE corner of outer D | 1996±16 |
| 3 times nearest distance between outer walls of Kivas B and C | 1929±15 |
| Weighted average | 1945±37 |

4. Discussion

For all measurements taken in the analysis, the ground survey confirmed the estimates made in the aerial survey within the statistical uncertainties, as seen in Table 1. The strong color gradations along the edges of site features in the aerial image of the Sun Temple were a significant aid in precise line placement, yielding statistical uncertainties on the measurements comparable to those achieved in the ground survey, with a relative uncertainty at the scale of the key site features of between 0.5% to 1%.

Aerial imagery is gaining increasing use in archaeology for location and analysis of sites (for instance, in References [28–31]). However, this is the first time, to our knowledge, that aerial imagery has been used to perform a survey of an archaeological site to the precision achieved here. The precision of the survey was aided by the good condition of the remaining site walls, and the availability of relatively high resolution aerial photographs of the site taken on clear days.

We found apparent evidence that a Golden rectangle was used to construct the outer D of the Sun Temple site, and that features of the geometrical construction of this rectangle are also associated with other geometrical constructs at the Sun Temple site. The Golden rectangle is evident in many examples of ancient Greek architecture [26]. However, evidence that societies in the prehistoric Americas knew of the Golden ratio is scant, although some evidence has been found that the footprints of examples of ancient Mayan ceremonial architecture may have exhibited approximations to the Golden rectangle [15].

We note here that, to within 1%, aerial survey measurements reveal that a Golden rectangle also encases the walls of another major ancestral Pueblo ceremonial site, Pueblo Bonito in Chaco Canyon, NM (see Figure 5). As shown in Figure 5, features related to the geometrical construction of that Golden rectangle are associated with other geometries of the site, like the arc of the northeast wall. Pueblo Bonito is much larger than the Sun Temple site (approximately 25 times larger in area), and was built over a period of several centuries from the 800's to around 1100 CE [32], ending around 100 to 150 years before Sun Temple was built. The final construction phases of Pueblo Bonito involved constructing the outermost walls that fit within the dimensions of a Golden rectangle, and the outer walls on the northeastern part of the site were shaped to conform to a circular arc associated with the Golden rectangle. Thus, it appears that the knowledge of how to construct Golden rectangles dated to at least 1100 CE in ancestral Puebloan society, and that the knowledge perhaps was shared through generations.

We also find evidence of at least one equilateral triangle (or, alternatively, a triangle with interior angles of 30°:60°:90°) is associated with the Sun Temple site, one side of which connects the Sun Shrine to the center of Kiva A, and another side of which is parallel to the south wall of the complex. Another equilateral triangle of height X/2 may be associated with the layout of the outer radius of Kiva D, inscribed within the triangle. A previous study of the Sun Temple site found that the small pecked basin to the north of the site formed a rough approximation to within about 5% to 10% of an equilateral triangle with the SW and SE corners of the D-shaped walls [12]. However, evidence has not been hitherto uncovered of societies in the prehistoric Americas employing more precise equilateral triangles (or 30°:60°:90° right triangles) in their architecture. The possible exception are the Maya, who are known to have constructed their round hearths with three stones placed approximately equidistant around the perimeter [33]; this arrangement was likely more practical than reflective of special attention to geometrical constructs, as even three roughly spaced stones help to prevent pots resting on the stones from tipping, similar to the benefits of a three-legged stool.

As shown in Figure 3, at least two Pythagorean 3:4:5 triangles appear to be associated with key features of the Sun Temple site. As shown in Figure 3, a third Pythagorean 3:4:5 triangle may be associated with Kiva D, the Sun Shrine, and the pecked basin to the north of the site, and the triangle is juxtapositionally associated with several other geometrical constructs at the site. For several millennia, the Pythagorean theorem has been known to scholars in China, India, Babylon and Greece [24]. There is also evidence that the ancient Egyptians in the time of the pharaohs were aware of at least certain cases of Pythagorean triples [34], but evidence of awareness of these triangles in the prehistorical Americas has hitherto not been uncovered.

We also find evidence that a common unit of measurement was used to layout many features of the Sun Temple site. The base unit is either $L = 30.5 \pm 0.5$ cm, or one third of that. Interestingly, several past societies in the world developed a unit of measure close to L, including the modern imperial foot, which is 30.48 cm, the Greek common foot of 31.50 cm, the Roman foot of 29.59 cm, and the "northern foot" of 33.53 cm (used particularly by Germanic peoples) [35, 36].

The average length of a modern male human foot is around 27 cm [37], which is less than these defined "foot" measurements, but is more consistent with the Pythic "natural foot" of 25 cm used by the ancient Celts [35]. It is hypothesized that the larger foot units were defined to more closely match integer multiples, or rational fractions, of other human body parts typically used for measurement, such as three times the width of a clenched fist (one "hand" is approximately 10 cm), or 2/3 times the length of a forearm from elbow to fingertip (one "cubit" is approximately 46 cm) [36]. Thus, in a similar vein, it may likely be that the base unit of measurement used to construct the Sun Temple was based on the width of a clenched fist, approximately 10 cm, with L being three times this measure. The prehistoric Maya had a unit of measure, *kab*, equal to 9.2 ± 0.3 cm, similar to the apparent ancestral Pueblo base unit of measure [38].

While there does appear to be a common unit of measurement to the layout of this particular site, it remains to be seen if this measure is evident at other prehistoric Pueblo sites of the same period. It may be that the unit of measurement used here was simply based on a body part of the site architect.

It is also interesting to note that we find multiples of 3, 4, and 12 of this base measurement in the site layout, similar to hypothesized base measurements underlying Mayan ceremonial architecture, which are also found in multiples of 3, 4, and 12 [38], despite the vigesimal number system of the Maya [39]; this perhaps reflects Mesoamerican influence on the Mesa Verde region [40, 41].

A previous study of a survey of the Sun Temple site by Munson *et al* (2010) mentioned the possibility that the site may have been laid out with some common unit of measurement (not specified in the study), perhaps using a base 10 mathematical system [12]. We find no evidence of multiples of 5 or 10 of the base unit identified by this study. This does not, of course, mean that multiples of 5 or 10 were not used at all in the site layout, merely that they were not apparent in the geometrical constructs identified in this study.

Given the previous observations of key solar and lunar alignments associated with site features [14], it is obvious that, similar to several examples of Mayan ceremonial architecture [42], and prehistoric mound builders in Ohio [43, 44], the planning of the Sun Temple site involved astronomy as well as attention to geometry. As Aveni (1982) pointed out when speaking of Mayan ceremonial architecture [42]:

*"Indeed, astronomical and geometrical determining factors, not all of them necessarily related, may have combined to influence the placement and orientation of a given architectural component."*

We note here that in the past, some individuals have pursued pseudo-scientific studies of geometric layouts of archaeological sites, and this unfortunately has led to a cachet of "woo science" associated with any study in this area, as is also unfortunately the case with the field of archaeoastronomy. As Professor Clive Ruggles of the University of Leicester once pointed out about the field of archaeoastronomy,

> "*A field with academic work of high quality at one end but uncontrolled speculation bordering on lunacy at the other.*"

The same is too often true of archaeogeometry; for instance, in various studies purporting evidence of geometric layouts associated with Southwest prehistoric sites (often associated with New Age theories), the author has noted that geometric shapes of arbitrary size are simply overlaid on a map of a site, without regard to whether or not the vertexes, sides, or size of the shape are meaningfully associated with any of the key features of the site [45], an example of "unanchored geometric interpretation" [46]. At times, many interconnecting lines associated with the vertexes of the geometric shapes are then also overlaid, making an image that looks impressively complex [47], but, in reality, is meaningless. In addition, at times rectangles of many various relative length to width ratios are used to tessellate features in a site map, under the premise that the tessellation is meaningful [48]; but, in fact, if you are given the option to choose from many different length to width ratios for the rectangles, you can tessellate nearly all rectangular sites with such shapes regardless whether the site was laid out with them or not. Further, the relative measurements of some purported geometries overlaid on crude site maps fail to be verified with examination using aerial imagery, for instance by using Google Earth. Most importantly, quantification of the statistical significance of the purported alignments or site measures is never provided in such studies.

To ensure rigor in our analysis, we thus only examined potential geometries associated with either the size of key features like the four Kivas and the outer D, or geometries associated with measures between at least two of the key features. We also only presented results that can be independently verified by interested readers using aerial imagery available with software programs such as Google Earth, and we assessed the statistical significance of the apparent geometries we examined. For simplicity, we constrained the analysis only to those elements that were found by previous analyses to have been built first (the four Kivas, and the outer walls of the outer D), and the Sun Shrine.

*Figure 5: Pueblo Bonito, a Chaco great house built by the ancestral Pueblo peoples between the 800's to early 1100's CE. Overlaid is a rectangle encasing the walls of the site. To within 1%, the rectangle is consistent with the dimensions of a Golden rectangle. The blue line shows the arc of a circle circumscribed within the square associated with the geometric construction of the rectangle. Google and the Google logo are registered trademarks of Google Inc., used with permission.*

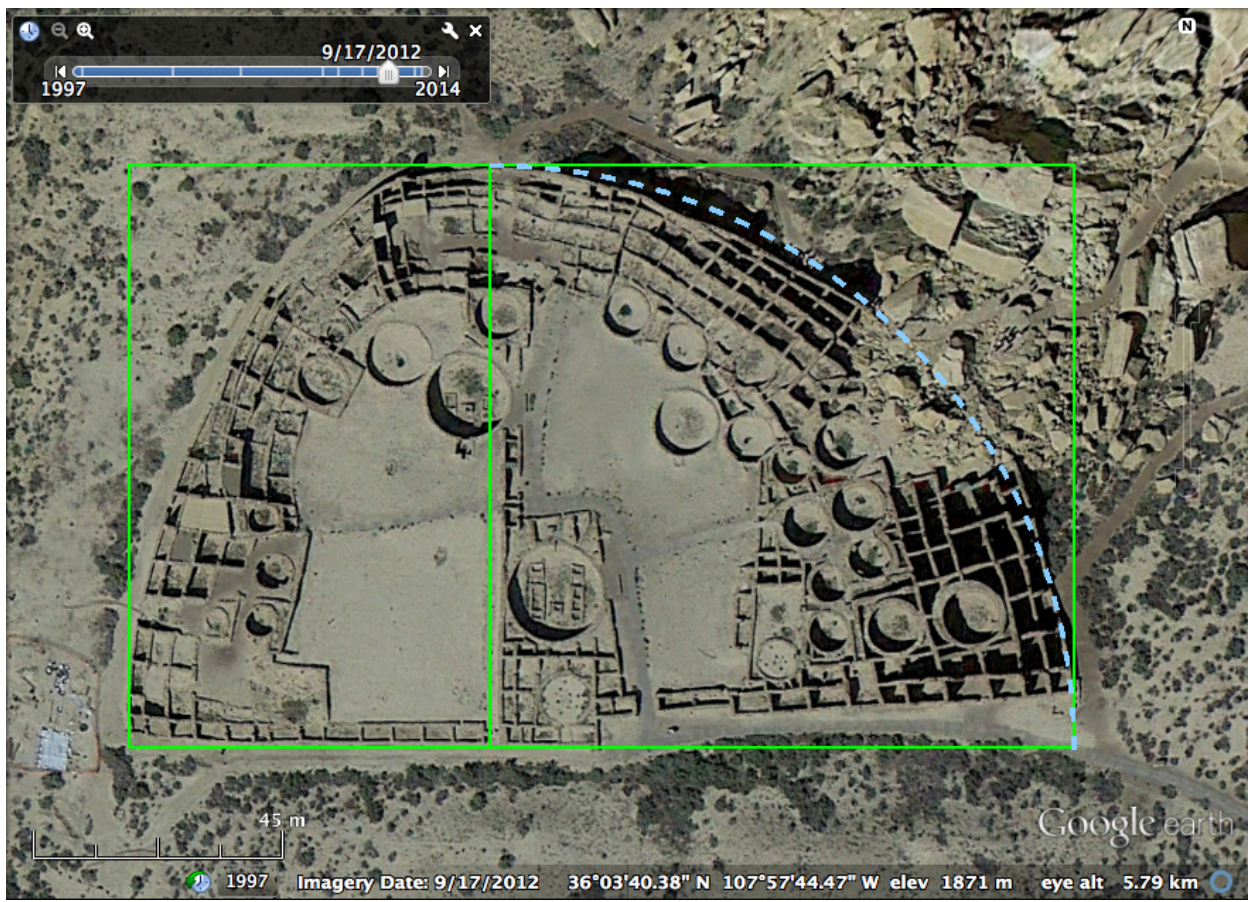

5. Summary

In this analysis, we examined the layout of the Sun Temple ceremonial complex in Mesa Verde National Park in Colorado, USA, through aerial and ground surveys. The site was built by the ancestral Pueblo peoples, c. 1200 CE.

We found evidence that key features of the site were apparently laid out using the Golden rectangle, squares, 45° triangles, Pythagorean 3:4:5 triangles, and equilateral triangles. A common unit of measurement appears to underlie these constructs, $L = 30.5 \pm 0.5$ cm, or an integer division of L. The site was laid out with remarkable precision, with the relative uncertainty on measurements estimated to be approximately 1%. Further study is needed to determine if this unit of measurement was particular to the Sun Temple (for instance, a dimension of a body part of the Sun Temple site architect), or whether the unit was common to other ancestral Pueblo sites. Further study is also needed to determine whether or not geometrical constructs are evident at other ancestral Pueblo sites other than Sun Temple and Pueblo Bonito.

By mere random chance, any site may yield potential evidence of a geometrical construct if enough site elements are examined, even when no such constructs were actually used in the design of the site. However, what makes this particular site unusual is the number of geometrical constructs found when just a few site elements were considered. And, most especially, the relationship of those geometric constructs to the apparent common unit of measurement at the site is extraordinarily unlikely to occur by mere random chance.

The Sun Temple site thus likely represents the first evidence of advanced knowledge of several geometrical constructs in prehistoric America. Given that the ancestral Pueblo peoples had no written language or number system, the precision of such a layout would be a remarkable feat. It is unclear why these ancients potentially felt the need to employ these constructs in the Sun Temple site. Perhaps the specialized knowledge of how to construct these shapes with a straightedge and a cord formed part of the inherent mysticism of the ceremonial nature of the site. Certainly, the care with which the site was laid out supports its role as a key center of ceremony and ritual in the region.

Accepted *Journal of Archaeological Science: Reports*, Jan 2017 (in press)

Acknowledgments

The author is grateful for the kind assistance of National Park Service, in particular NPS archaeologist Gay Ives, during the site work related to this study, and to Greg Munson for very helpful discussions.

References

[1] Gilbert R Wenger. *The Story of Mesa Verde National Park*. Mesa Verde Museum Association, 1987.

[2] Mark A Stiger. Mesa Verde subsistence patterns from Basketmaker to Pueblo III. *The Kiva*, pages 133–144, 1979.

[3] Michael A Adler, Todd Van Pool, and Robert D Leonard. Ancestral Pueblo population aggregation and abandonment in the North American southwest. *Journal of World Prehistory*, 10(3):375–438, 1996.

[4] Linda S Cordell, Carla R Van West, Jeffrey S Dean, and Deborah A Muenchrath. Mesa Verde settlement history and relocation: Climate change, social networks, and ancestral Pueblo migration. *The Kiva*, pages 379–405, 2007.

[5] Mark Varien. *Sedentism and Mobility in a Social Landscape: Mesa Verde & Beyond*. University of Arizona Press, 1999.

[6] Carla R Van West and Jeffrey S Dean. Environmental characteristics of the CE 900-1300 period in the central Mesa Verde region. *The Kiva*, pages 19–44, 2000.

[7] Donna M Glowacki. *Living and Leaving: A Social History of Regional Depopulation in Thirteenth-Century Mesa Verde*. University of Arizona Press, 2015.

[8] Harry T Getty. New dates from Mesa Verde. *Tree-Ring Bulletin*, 1(3):21–23, 1935.

[9] Guy E Gibbon and Kenneth M Ames. *Archaeology of prehistoric native America: an encyclopedia*, volume 1537. Taylor & Francis, 1998.

[10] Jesse Walter Fewkes. *Excavation and repair of Sun Temple, Mesa Verde National Park*. Govt. print. off., 1916.

[11] Gregory E Munson. Legacy documentation: using historical resources in a cultural astronomy project. *Proceedings of the International Astronomical Union*, 7(S278): 265–274, 2011.

[12] Gregory Munson, Larry Nordby, and Bryan Bates. Reading, writing, and recording the architecture: How astronomical cycles may be reflected in the architectural construction at Mesa Verde National Park, Colorado. *Archeoastronomy*, 23:44, 2010.

[13] National Park Service. *Mesa Verde National Park, Colorado*. U.S. Department of the Interior, 2011.

[14] J McKim Malville and Claudia Putnam. *Prehistoric astronomy in the Southwest*. Big Earth Publishing, 1993.

[15] James A Doyle. Early Maya geometric planning conventions at El Palmar, Guatemala. *Journal of Archaeological Science*, 40(2):793–798, 2013.

[16] Francine Vinette. In search of Mesoamerican geometry. *Native American mathematics*, page 387, 1986.

[17] Francine Vinette. *Manifestations of geometrical concepts in Mesoamerica*. PhD thesis, University of Ottawa, 1982.

[18] Jason Ur. Google Earth and archaeology. *The SAA Archaeological Record*, 6(3):35–38, 2006.

[19] Chen Zhao. Introduction of Google Earth's practical use in archaeology. *Southeast Culture*, 2:006, 2007.

[20] Adrian Myers. Field work in the age of digital reproduction: a review of the potentials and limitations of Google Earth for archaeologists. *SAA Archaeological Record*, 10(4):7–11, 2010.

[21] Karim Sadr and Xavier Rodier. Google Earth, GIS and stone-walled structures in southern Gauteng, South Africa. *Journal of Archaeological Science*, 39(4): 1034–1042, 2012.

[22] David Kennedy and MC Bishop. Google Earth and the archaeology of Saudi Arabia. a case study from the Jeddah area. *Journal of Archaeological Science*, 38(6): 1284–1293, 2011.

[23] Michael L Hargrave. Ground truthing the results of geophysical surveys. *Remote Sensing in Archaeology: An Explicitly North American Perspective*, pages 269–304, 2006.

[24] Bartel L Van der Waerden. *Geometry and algebra in ancient civilizations*. Springer Science & Business Media, 2012.

[25] George E Martin. *Geometric constructions*. Springer Science & Business Media, 2012.

[26] Mario Livio. *The golden ratio: The story of phi, the world's most astonishing number*. Broadway Books, 2008.

[27] Olive Jean Dunn. Multiple comparisons among means. *Journal of the American Statistical Association*, 56(293):52–64, 1961.

[28] Jesse Casana and Jackson Cothren. Stereo analysis, DEM extraction and orthorectification of CORONA satellite imagery: archaeological applications from the Near East. *Antiquity*, 82(317): 732–749, 2008.

[29] Robert H Bewley, SP Crutchley, and CA Shell. New light on an ancient landscape: Lidar survey in the Stonehenge World Heritage Site. *Antiquity*, 79(305): 636–647, 2005.

[30] Véronique De Laet, Etienne Paulissen, and Marc Waelkens. Methods for the extraction of archaeological features from very high-resolution Ikonos-2 remote sensing imagery, Hisar (southwest Turkey). *Journal of Archaeological Science*, 34(5): 830–841, 2007.

[31] Mark Altaweel. The use of ASTER satellite imagery in archaeological contexts. *Archaeological Prospection*, 12(3): 151–166, 2005.

[32] Thomas C Windes and Dabney Ford. The Chaco wood project: The chronometric reappraisal of Pueblo Bonito. *American Antiquity*, pages 295–310, 1996.

[33] Dorie Reents-Budet. Elite Maya pottery and artisans as social indicators. *Archeological Papers of the American Anthropological Association*, 8(1): 71–89, 1998.

[34] Beatrice Lumpkin. The Egyptians and Pythagorean triples. *Historia Mathematica*, 7(2): 186–187, 1980.

[35] June A Sheppard. Metrological analysis of regular village plans in Yorkshire. *The Agricultural History Review*, 22(2): 118–135, 1974.

[36] Herbert Arthur Klein. *The science of measurement: A historical survey*. Courier Corporation, 2012.

[37] Claire C Gordon, Thomas Churchill, Charles E Clauser, Bruce Bradtmiller, and John T McConville. Anthropometric survey of US army personnel: methods and summary statistics 1988. Technical report, DTIC Document, 1989.

[38] Patricia J O'Brien and Hanne D Christiansen. An ancient Maya measurement system. *American antiquity*, pages 136–151, 1986.

[39] Gary D Salyers. The number system of the Mayas. *Mathematics Magazine*, 28(1): 44–48, 1954.

[40] Carroll L Riley. *Becoming Aztlan: Mesoamerican influence in the greater Southwest, CE 1200-1500*. Univ of Utah Pr, 2005.

[41] Randall H McGuire. The Mesoamerican connection in the southwest. *The Kiva*, pages 3–38, 1980.

[42] Anthony F Aveni and Horst Hartung. Precision in the layout of Maya architecture. *Annals of the New York Academy of Sciences*, 385(1): 63–80, 1982.

[43] Ray Hively and Robert Horn. Geometry and astronomy in prehistoric Ohio. *Journal for the History of Astronomy Supplement*, 13:1, 1982.

[44] Ray Hively and Robert Horn. Hopewellian geometry and astronomy at High Bank. *Journal for the History of Astronomy Supplement*, 15:85, 1984.

[45] Gary David. *The kivas of heaven*. Adventures Unlimited Press, 2010.

[46] Stephen Skinner. *Sacred geometry: deciphering the code*. Sterling Publishing Company, Inc., 2009.

[47] Charly Gullett. *Engineering Chaco, the Anasazi field manual*. Warfield Press, 2011.

[48] Christopher Powell. *The shapes of sacred space: a proposed system of geometry used to lay out and design Maya art and architecture and some implications concerning Maya cosmology*. PhD thesis, University of Texas at Austin, 2010.